\documentclass[12pt,a4paper]{article}

\usepackage[headinclude]{typearea}
\usepackage[english]{babel}           
\usepackage[utf8]{inputenc}
\usepackage[T1]{fontenc}
\usepackage{scrtime}
\usepackage{srcltx}
\usepackage{amsmath,amsthm,amsfonts,amssymb}
\usepackage {color}                                  
\usepackage {pifont}                                 
\usepackage {multicol,multirow}                               
\usepackage[numbers]{natbib}                                                     
\usepackage[normalem]{ulem}                                                     
\usepackage [scaled=.90]{helvet}                     
\usepackage {courier}                                
\usepackage {graphicx}                               
\usepackage {array}                                  
\usepackage {longtable}                              
\usepackage{fancyvrb}
\usepackage{enumerate} 
\usepackage{ifpdf}
\usepackage{caption}
\usepackage{mathrsfs}

\usepackage{setspace}

\usepackage{marvosym}

\theoremstyle{plain}
\newtheorem{theorem}{Theorem}[section]

\newtheorem{lemma}[theorem]{Lemma}

\newtheorem{proposition}[theorem]{Proposition}

\theoremstyle{definition}

\newtheorem{definition}[theorem]{Definition}

\newtheorem{algorithm}[theorem]{Algorithm}

\newcommand{\exclude}[1]{}

\newcommand{\othereps}{\zeta}
\newcommand{\1}{{\mathbf{1}}}
\newcommand{\E}{{\mathbb{E}}}

\newcommand{\N}{{\mathbb{N}}}
\renewcommand{\P}{{\mathbb{P}}}

\newcommand{\R}{{\mathbb{R}}}

\definecolor{darkgreen}{rgb}{0,0.5,0}
\definecolor{lightgreen}{rgb}{0.5,0.9,0.5}
\definecolor{magenta}{rgb}{0.75,0,0.25}

\definecolor{violet}{rgb}{0.25,0,0.75}

\newcommand{\hypsurf}{\Theta}
\newcommand{\appx}{X^\delta}
\newcommand{\appz}{Z^\delta}
\newcommand{\es}{{\underline{s}}}

\renewcommand{\P}{{\mathbb P}}

\newcommand{\cS}{{\cal S}}

\newcommand{\be}{\begin{equation}}
\newcommand{\ee}{\end{equation}}
\newcommand{\bea}{\begin{eqnarray}}
\newcommand{\eea}{\end{eqnarray}}
\newcommand{\beast}{\begin{eqnarray*}}
\newcommand{\eeast}{\end{eqnarray*}}
\newcommand{\bproof}{\begin{proof}}
\newcommand{\eproof}{\end{proof}}

\begin{document}

\title{Numerical methods for SDEs with drift  discontinuous
on a set of positive reach}

\author{Gunther Leobacher \and Michaela Sz\"olgyenyi}

\date{Corrected version, December 2018}

\maketitle

\begin{abstract}
For time-homogeneous stochastic differential equations (SDEs) it is enough 
to know that the coefficients are Lipschitz to conclude existence and 
uniqueness of a solution, as well as the existence of a strongly convergent
numerical method for its approximation.
Here we introduce a notion of piecewise Lipschitz functions and study SDEs with a drift coefficient satisfying only
this weaker regularity condition.
For these SDEs we can construct a strongly convergent approximation scheme,
if the set of discontinuities is a sufficiently smooth hypersurface satisfying the geometrical property of being
of positive reach.
We then arrive at similar conclusions as in the Lipschitz case.
We will see that, although SDEs are in the center of our interest, we  will
talk surprisingly little about probability theory here.
\end{abstract}

\paragraph{About the authors:} \emph{Gunther Leobacher} studied mathematics at the University of Salzburg,
and finished his PhD in 2001 under the supervision of G.~Larcher at Johannes Kepler University Linz.
In 2002--2003 he was a Postdoc with L.G.C.~Rogers at the University of Cambridge. In 2012 he became an as\-so\-ci\-ate professor at JKU Linz.
In February 2017 he was appointed full professor
of stochastics at the University of Graz.

\emph{Michaela Sz\"olgyenyi} studied mathematics at the Johannes Kepler University Linz, and finished her PhD
in 2015 under the supervision of G.~Leobacher. Then she became a Postdoc with R.~Frey
at Vienna University of Economics and Business. From August 2017 she works at ETH Z\"urich in the research group of A.~Jentzen.
In 2017 she was granted an international research project by the AXA Research Fund.\\

\section{Introduction}

Stochastic differential equations (SDEs) are essential for many models in mathematical finance, risk theory, biology, physics, and chemistry. Usually, these equations cannot be solved explicitly.
Hence, we are interested in finding numerical methods with positive convergence speed for solving them.

We consider general SDEs on the $\R^d$, which are of the form
\begin{align}\label{eq:SDE}
d X_t = \mu(X_t) dt +\sigma(X_t) dW_t\qquad
X_0=x\,,
\end{align}
with initial value $x\in \R^d$,
drift coefficient $\mu:\R^d\longrightarrow \R^d$,
diffusion coefficient $\sigma:\R^d\longrightarrow \R^{d\times m}$, and
$m$-dimensional standard Brownian motion $W$ (thus adding noise to the ordinary differential equation).
Little generality is lost if we assume $m=d$, and we will do so throughout this article.

By a (strong) solution we mean a continuous stochastic process $X$ that
is adapted to the filtration generated by $W$ and that satisfies
\begin{align}\label{eq:SDEint}
X_t=x+\int_0^t \mu(X_s) ds + \int_0^t \sigma(X_s)dW_s
\end{align}
for all $t\ge 0$ almost surely. The solution $X$ is unique, if the paths of any
other solution to \eqref{eq:SDEint} coincide with those of $X$ almost surely.

The second integral in \eqref{eq:SDEint} is It\^o's stochastic integral, the construction of
which we will not repeat here. 	
Suffice it to mention that for $K$ from a suitable class of stochastic processes it holds that
\[
\int_0^t K_s \, dW_s=\lim_{n\to \infty}\sum_{k=1}^{\lfloor 2^n t\rfloor-1}
K_{k 2^{-n}} (W_{(k+1) 2^{-n}}-W_{k 2^{-n}})\,,
\]
reminding us of the Riemann integral (but with evaluation of the integrand only in the left boundary of small intervals).
A particularity of It\^o's integral is that there appears a correction term in the fundamental theorem of calculus, that is, for $X_t=X_0 + \int_0^t H_s ds + \int_0^t K_s dWs$ and for a sufficiently regular function $f:\R \longrightarrow \R$,
\begin{align*}
 f(X_t)=f(X_0)+\int_0^t f'(X_s) H_s ds +\int_0^t f'(X_s) K_s dW_s + \frac{1}{2} \int_0^t f''(X_s) K_s^2 ds\,.
\end{align*}
This is known as It\^o's formula.
The rigorous construction of the stochastic integral gave
meaning to the concept of a solution of an SDE. In addition to that \citet{ito1951} proved that
a unique solution to \eqref{eq:SDE} exists, whenever $\mu$ and $\sigma$ are Lipschitz-continuous.

Under the same assumptions \citet{maruyama1955} proved that the Euler-Maruyama (EM) scheme
\begin{align*}
\appx_t=x+\int_0^t \mu(\appx_\es) ds + \int_0^t \sigma(\appx_\es) dW_s\,,
\end{align*}
with $\es=j\delta$ for $s\in[j \delta , (j+1)\delta )$, $j=0,\dots,(T-\delta)/\delta$,
(which reminds us of the Euler scheme for ordinary differential equations, but with an additional term corresponding to the stochastic integral)
converges with strong order $1/2$. In general we say that 
a numerical approximation $\appx$ converges with strong order $\gamma$, if for any fixed $T>0$, there exists a constant $C$ such that for 
sufficiently small step-size $\delta>0$ 
it holds that   
\[
\E\Big(\sup_{0\le t \le T}\|X_t-\appx_t\|^2\Big)^{1/2}\le C \delta^{\gamma}\,.
\]
Higher order algorithms exist under stronger regularity conditions on the coefficients,
most notably the Milstein method and stochastic Runge-Kutta schemes, see \citet{kloeden1992}.

The question of how to solve SDEs with irregular (non-globally Lipschitz) coefficients approximately is a very active topic of research.
There is still a big gap between the assumptions on the coefficients of
these equations under which strong convergence with convergence rate has been
proven in the scientific literature, and the assumptions that equations in
real-world applications satisfy.

In contrast to that, several delimiting results have been proven recently, stating that a certain SDE with relatively well-behaved (infinitely often differentiable) coefficients cannot be solved approximately in finite time, cf.~\citet{hairer2015,jentzen2016,muellergronbach2016,yaroslavtseva2016}.
However, there is still a big discrepancy between the assumptions on the coefficients under which convergence with strong convergence rate has been proven and the properties of the coefficients of the SDE presented in \citet{hairer2015}.

Here we narrow the gap described above by settling convergence with positive convergence speed of a numerical method for $d$-dimensional SDEs with discontinuous drift and degenerate diffusion coefficient.
First steps in this direction have previously been made by
\citet{ngo2016c}, who proved convergence of order up to $1/4$ of the Euler-Maruyama method for $d$-dimensional SDEs which have a discontinuous, bounded drift that satisfies a one-sided Lipschitz condition and a H\"older continuous, bounded, and uniformly non-degenerate diffusion coefficient.
In \citet{ngo2016a,ngo2016b} they do not need the one-sided Lipschitz condition any more, but the result only works for one-dimensional SDEs and relies on uniform non-degeneracy of the diffusion coefficient.

SDEs with discontinuous drift appear naturally when studying stochastic optimal control problems with bang-bang type optimal strategies,
that is with strategies of the form $\1_\cS(X)$ for a measurable set $\cS \subseteq \R^d$.
If in addition only a noisy signal of the underlying state process $X$ is available, then filtering this signal leads to a degenerate diffusion coefficient and increases the dimension substantially.
Examples can be found in \citet{sass2004,rieder2005,frey2012,sz12,sz16,sz2016a,shardin2017}.

\begin{figure}[ht!]
 \begin{center}
\scalebox{1.0}{\input{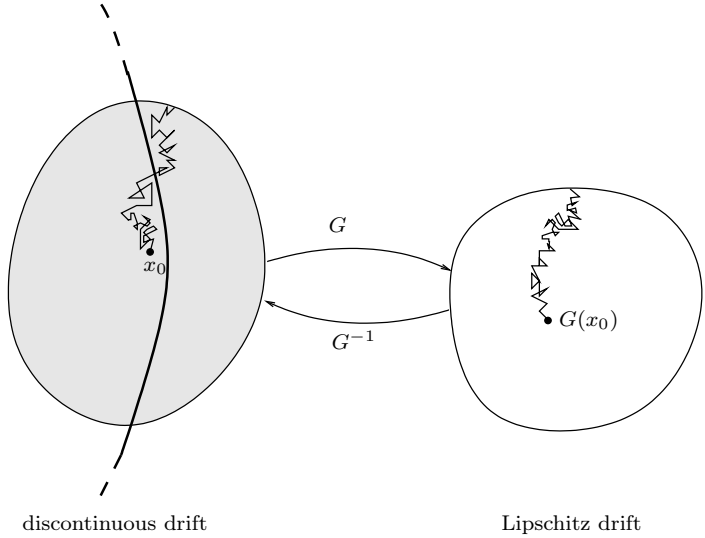}}
\caption{A sketch of the idea for the construction of our numerical method.\label{fig:idea}}
\end{center}
\end{figure}

The idea for tackling the problem is illustrated in Figure \ref{fig:idea}:
to overcome the issues caused by a discontinuous drift coefficient, we want to find a transform $G$ with the property that the coefficients of the transformed SDE for $G(X)$ are Lipschitz.
Then we want to apply the EM scheme to that SDE, which converges with strong order $1/2$, to obtain an approximation of the solution to the transformed SDE.
In the end, we want to transform back to obtain an approximation of the solution to the original SDE \eqref{eq:SDE}.
In Figure \ref{fig:idea}, the set of discontinuities of the drift is illustrated by a smooth curve. Indeed, we need to make some assumptions to that end 
so that we can carry through our idea.

Thus, we have to solve the following tasks: 

\begin{enumerate}\itemsep-\parsep
 \item construct $G$ and all prove necessary properties;
 \item prove an existence and uniqueness result;
 \item construct a numerical method using $G$ (called GM) and prove convergence and convergence rate;
 \item prove convergence and convergence rate for EM starting from GM.
\end{enumerate}

We will start with presenting our results in dimension one, and subsequently we will show how these ideas can be extended to general dimension.

This is a review article; the results and examples presented here, can 
be found in \citet{sz15,sz2016b,sz2017c}.

\section{Result in dimension one}

In order to construct an appropriate $G:\R \longrightarrow \R$, 
we have to know how such a transform
acts on the coefficients:  assuming existence of a solution $X$ and also
validity of It\^o's formula for $X$ and $G$ 
we get
\begin{align*}
G(X_t)=G(x)+\!\int_0^t G'(X_s)\mu(X_s) ds + \!\int_0^t G'(X_s)\sigma(X_s)dW_s
+\!\int_0^t \frac{1}{2}G''(X_s)\sigma(X_s)^2 ds\,.
\end{align*}
Thus $Z=G(X)$ is the solution of an SDE
with coefficients 
\begin{align*}
\tilde \mu (Z) &= G'(G^{-1}(Z))\mu(G^{-1}(Z))+\frac{1}{2}G''(G^{-1}(Z))\sigma(G^{-1}(Z))^2\,,\\
\tilde \sigma (Z) &= G'(G^{-1}(Z))\sigma(G^{-1}(Z))\,.
\end{align*}
Hence, $G$ maps  $X \mapsto Z$ and it
transforms $\mu,\sigma$ into $\tilde \mu,\tilde \sigma$.

We see that if $G\in C^2$ -- the classical assumption for It\^o's formula --
then $\tilde \mu,\tilde \sigma$
are continuous, if and only if $\mu,\sigma$ are continuous.
However 
if $G\in C^1$ and $\sigma$ is continuous and non-zero,
then we can offset jumps of $\mu$  with jumps of $G''$.
So we can get continuous $\tilde\mu$ from discontinuous $\mu$ with a less smooth transform.
Note that $\tilde \sigma$ is continuous in either case.
Hence, we choose $G\in C^1$ to be able to eliminate the discontinuities from the drift.
Note that we will have to verify that the heuristic application of It\^o's formula above is valid,
since the classical It\^o formula holds for $C^2$ functions.

With this, we are able to relax the Lipschitz condition on the drift.
\begin{definition}\label{def:pwlip1}
A function $\mu:\R\longrightarrow \R$ is called {\em piecewise Lipschitz}, if
there are finitely many points $\xi_1<\dots<\xi_m$  such that 
the restriction of $\mu$ to each of the intervals 
$(-\infty,\xi_1),(\xi_m,\infty)$ and $(\xi_k,\xi_{k+1})$, $k=1,\dots,m-1$,
is Lipschitz.
\end{definition}

For the presentation here, we now assume that $\mu$ is piecewise Lipschitz with only one jump in $\xi$,
but note that our result also holds for multiple jumps. Let
\begin{itemize}\setlength{\itemsep}{0em}
 \item $\mu$ be Lipschitz on $(-\infty,\xi)$ and $(\xi,\infty)$;
 \item $\sigma:\R\longrightarrow\R$ be Lipschitz with $\sigma(\xi)\neq 0$.
\end{itemize}
Note that the last condition is by far weaker than uniform non-degeneracy,
as for non-degeneracy one would need $\sigma$ to be bounded away from $0$ on the whole of $\R$.

We define the transform $G:\R\longrightarrow \R$ by
\begin{align}
\label{eq:transf1d}G(x)=x+\alpha (x-\xi)|x-\xi| \phi\left(\frac{x-\xi}{c}\right)=:x+\alpha \bar \phi(x)\,,
\end{align}
where $\alpha,c$ are appropriate constants, and
\[
\phi(u)=
\begin{cases}
(1+u)^3(1-u)^3 & \text{if } |u|\le 1\,,\\
0 & \text{else}
\end{cases}
\]
localizes the impact of $G$.
If $0<c<1/6|\alpha|$, then $G'>0$, and hence $G$ is globally invertible. Furthermore, we can prove that $G$ and $G^{-1}$ are Lipschitz.

Setting $Z=G(X)$, we have
\begin{align*}
dZ_t=\tilde \mu(Z_t) dt +\tilde\sigma(Z_t) dW_t\,,
\end{align*}
where
\begin{align*}
\tilde\mu(z)&=\mu(G^{-1}(z))+\frac{1}{2}\alpha\bar\phi''(G^{-1}(z))\sigma(G^{-1}(z))^2+\alpha\bar \phi'(G^{-1}(z)) \mu(G^{-1}(z))\,,\\
\tilde\sigma(z)&=\sigma(G^{-1}(z)) +\alpha\bar \phi'(G^{-1}(z))\sigma(G^{-1}(z))\,.
\end{align*}

In order to offset the jump of $\mu$ in $\xi$ by the jump of $G''$ (by construction also in $\xi$), we choose $\alpha$ as
\[
\mu(\xi+)+\frac{1}{2}\alpha\bar\phi''(\xi+)\sigma(\xi)^2
=\mu(\xi-)+\frac{1}{2}\alpha\bar\phi''(\xi-)\sigma(\xi)^2 \Longrightarrow
\alpha=\frac{\mu(\xi-)-\mu(\xi+)}{2\sigma(\xi)^2}\,.
\]
With this choice of $\alpha$ we have that $\tilde \mu$ is continuous.
\begin{lemma}[Elementary but essential]\label{lem:elem-1}
Let $\tilde \mu:\R\longrightarrow\R$ be a function satisfying
\begin{enumerate} \itemsep-\parsep
\item $\tilde \mu$ is continuous;
\item $\tilde \mu$ is piecewise Lipschitz.
\end{enumerate}
Then $\tilde \mu$ is Lipschitz.
\end{lemma}
Altogether we have that the coefficients of the SDE for $Z$ are Lipschitz.

Now, we are ready to prove the following theorem.
\begin{theorem}[\citet{sz15}]
 Let $\mu$ be piecewise Lipschitz and let $\sigma$ be Lipschitz and $\mu(\xi+)\ne\mu(\xi-)\Longrightarrow \sigma(\xi)\ne 0$.
 
 Then there exists a unique strong solution to the one-dimensional version of
\eqref{eq:SDE}.
\end{theorem}

The proof works as follows:
\begin{itemize}\itemsep-\parsep
 \item show that the SDE for $Z=G(X)$ has Lipschitz coefficients using Lemma \ref{lem:elem-1};
 \item then by It\^o's theorem, there exists a unique strong solution to this SDE; 
 \item set $X=G^{-1}(Z)$ and apply It\^o's formula to it, to see that
\begin{align*}
d X_t = \mu(X_t) dt +\sigma(X_t) dW_t\,.
\end{align*}
\end{itemize}
So we have constructed a process $X$ that solves our SDE.
There is one issue that we have already mentioned above: $G^{-1}\notin C^2$.
But in 1D, It\^o's formula holds nevertheless, see \cite[Problem 7.3]{karatzas1991}.

As sketched in Figure \ref{fig:idea} above, the transformation method in a natural way also leads to the following numerical scheme.

\goodbreak

\begin{algorithm}[\citet{sz15}]\label{alg:GM}
Given $\mu,\sigma, x,T$, and the step-size $\delta>0$,
\begin{enumerate}\itemsep-\parsep
\item precompute $G,G^{-1},\tilde \mu, \tilde \sigma$;
\item solve $dZ=\tilde \mu(Z)dt+\tilde \sigma(Z) dW$, $Z_0=G(x)$ on $[0,T]$
using the EM method to obtain the EM approximation
$\appz$;
\item compute the numerical approximation 
$\bar X_t=G^{-1}(\appz_t)$, for $t\in [0,T]$.
\end{enumerate}
\end{algorithm}

\begin{theorem}[\citet{sz15}]
 Let $\mu$ be piecewise Lipschitz and let $\sigma$ be Lipschitz and $\mu(\xi+)\ne\mu(\xi-)\Longrightarrow \sigma(\xi)\ne 0$.
 
 Then Algorithm \ref{alg:GM} converges with strong order $1/2$.
\end{theorem}

The proof is straightforward:
\citet{maruyama1955} showed that for sufficiently small step-size $\delta>0$,
\[
\E\Big(\sup_{0\le t\le T}|Z_t-\appz_t|^2\Big)^{1/2}\le C \delta^{1/2}\,.
\]

Denote by $L_{G^{-1}}$ the Lipschitz constant of $G^{-1}$. We get
\begin{align*}
\E\Big( \sup_{0\le t\le T}|X_t-\bar X_t|^2\Big)^{1/2}
&=\E\Big( \sup_{0\le t\le T}|G^{-1}(Z_t)-G^{-1}(\appz_t)|^2\Big)^{1/2}\\
&\le {L_{G^{-1}}}\, \E\Big( |Z_t-\appz_t|^2\Big)^{1/2} 
\le {L_{G^{-1}}} \, C \delta^{1/2}\,. 
\end{align*}

The following theorem is of particular relevance for the practical
implementation and efficiency of our
algorithm and shows that in 1D our result is already quite satisfactory:

\begin{theorem}[\citet{sz15}]\label{thm:inverseexp}
We can define an alternative transform $\widehat G$ which fulfills all the
necessary properties and which is piecewise cubic.
\end{theorem}

The relevance of this theorem lies in the fact that for a piecewise
cubic function the inverse can easily be computed explicitly.

\section{Result in general dimension}
\label{sec:ddim}

Extending our results to the multidimensional setting poses several challenges:
\begin{enumerate}\itemsep-\parsep
\item introduce a notion of \emph{piecewise Lipschitz} functions;
\item prove that piecewise Lipschitz $+$ continuous implies Lipschitz;
\item find a transform $G$ that makes the drift continuous;
\item show that $G$ has a global inverse;
\item show that It\^o's formula holds for $G^{-1}$.
\end{enumerate}

In this section we will sketch how these challenges were addressed.

\subsection{Piecewise Lipschitz functions on the $\R^d$}

There is no unique or universally accepted notion of a piecewise Lipschitz 
function on a subset of the $\R^d$. Below we propose such a definition that
generalizes the one-dimensional notion.

We call a continuous function $\gamma:[0,1]\longrightarrow A\subseteq \R^d$ a {\em curve 
in $A$ from $\gamma(0)$ to $\gamma(1)$} and we denote by
\[\ell(\gamma):=\sup\Big\{\sum_{k=1}^n |\gamma(t_{k})-\gamma(t_{k-1})|:
n\in\N,\,0=t_0<\dots<t_n=1\Big\}\]
its (possibly infinite) length.

\begin{definition}\label{def:intrinsic}
Let $\emptyset \ne A\subseteq \R^d$. Define the {\em intrinsic metric} on
$A$ by 
\[
\rho(x,y):=\inf\{\ell(\gamma): \text{$\gamma$ a curve in $A$ from $x$ to $y$} \}\,.
\]
Here, the infimum over an empty set is defined as $\infty$. 
\end{definition}

\begin{definition}\label{def:pwlip2}
A function $\mu:\R^d\longrightarrow \R^m$ is {\em piecewise Lipschitz},
 if there exists a hypersurface $\Theta$ with finitely 
many connected components such that the restriction $\mu|_{\R^d\backslash \Theta}$ 
is Lipschitz w.r.t.~the 
intrinsic metric on $\R^d\backslash \Theta$, and w.r.t.~the Euclidean metric on $\R^m$.

In that case we call $\Theta$ an {\em exceptional set} for $\mu$.
\end{definition}

Note that the definition coincides with Definition  
\ref{def:pwlip1} for $d=1$. It shares also some basic and well-known properties
with the elementary definition.

\begin{proposition}\label{prop:difflip}
Let $\Theta$ be a hypersurface in $\R^d$ and let 
$\mu:\R^d\longrightarrow \R^m$ be a function such that 
$\mu|_{\R^d\backslash\Theta}$ is differentiable with bounded derivative.

Then $\mu$ is piecewise Lipschitz with exceptional set $\Theta$
and 
\[\sup_{x,y\in \R^d\backslash \Theta:\rho(x,y)>0}\frac{\|\mu(x)-\mu(y)\|}{\rho(x,y)}=\sup_{x\in \R^d\backslash \Theta}\|\mu'(x)\|\,.\] 
\end{proposition}

The following lemma is almost trivial in dimension one, but not so
in general dimension: 

\begin{lemma}\label{lem:elem-multi}
Let $\mu:\R^d\longrightarrow \R^m$ be a function such that
\begin{enumerate}\itemsep-\parsep
\item $\mu$ is continuous;
\item $\mu$ is piecewise Lipschitz with exceptional set $\Theta$;
\item\label{it:finite} $\Theta$ is such that  for all
$x,y\in \R^d\backslash \Theta$ and  all $\eta>0$ there exists a 
curve $\gamma$ in the $\R^d$ from $x$ to $y$
such that $\ell(\gamma)<\|y-x\|+\eta$ and
$\#\left(\gamma\cap\Theta\right)<\infty$.
\end{enumerate}
Then $\mu$ is Lipschitz (w.r.t. the Euclidean norm) with Lipschitz constant
\[
L_\mu=\sup_{x,y\in \R^d\backslash \Theta:\rho(x,y)>0}\frac{\|\mu(x)-\mu(y)\|}{\rho(x,y)}\,.
\]
\end{lemma}

Lemma \ref{lem:elem-multi} differs from Lemma \ref{lem:elem-1} essentially 
by item \ref{it:finite}, which is trivially satisfied in dimension one by our 
definition of `piecewise Lipschitz'.
\begin{figure}
\begin{center}
\includegraphics[scale=0.5]{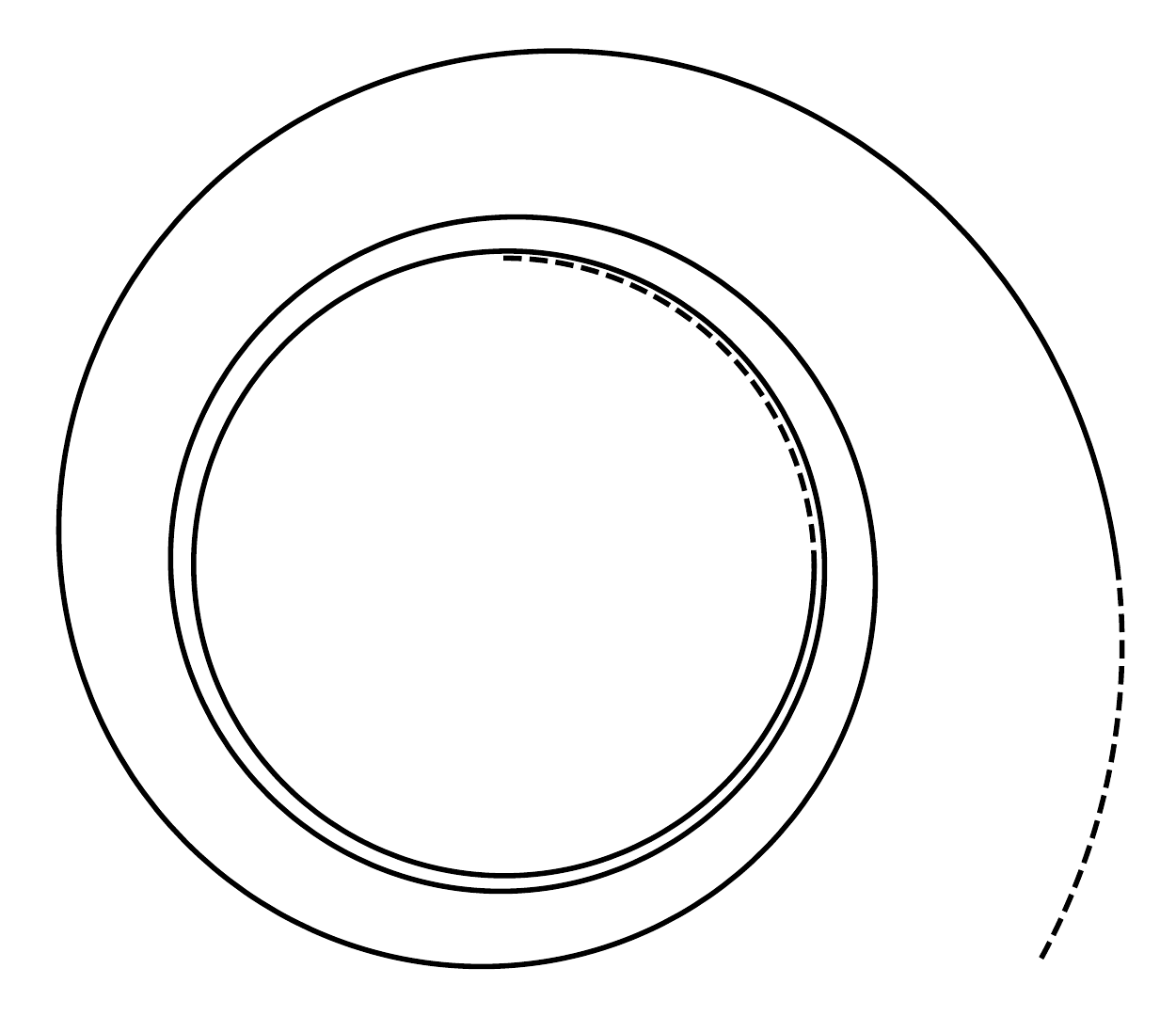}
\end{center}
\caption{An example for a hypersurface with bounded derivative of the 
unit normal vector which is not of positive reach.}\label{fig:bounded-ndash}
\end{figure}

Figure \ref{fig:bounded-ndash} shows an example of a two-dimensional $C^\infty$-hypersurface (i.e.~a curve) for which item \ref{it:finite} of 
Lemma \ref{lem:elem-multi} is not satisfied. The following notion
will prove useful for this issue:

\begin{definition}
A subset $\Theta\subseteq\R^d$ is of 
\emph{positive reach}, if there exists $\varepsilon>0$ such that for every $x\in
\R^d$ with $d(x,\Theta)<\varepsilon$ there is a unique $p\in \Theta$ with
$\|x-p\|=d(x,\Theta):=\inf\{\|x-\xi\|:\xi \in \Theta\}$. 
\end{definition}

If $\Theta$ has positive reach, then the projection map $p$ which assigns to 
$x$ its closest point $p(x)$ on $\Theta$ is a well-defined single-valued map 
on $\Theta^\varepsilon:=\{x\in \R^d:d(x,\Theta)<\varepsilon\}$ for some $\varepsilon>0$. 
Examples of
hypersurfaces having this property include hyperplanes and all compact
$C^2$-hypersurfaces,
which follows from the lemma in \citet{foote1984}, where it is also shown 
that the projection map $p$ is in $C^{k-1}$ if $\Theta$ is in $C^k$.
We will always assume that the set of discontinuities
of the drift coefficient  is of positive reach. 

The projection map $p$ will play a prominent role in the construction of 
the multivariate transform $G$.

One consequence of the positive reach property for $\Theta$ is that item 
\ref{it:finite} of Lemma \ref{lem:elem-multi} is automatically satisfied.
This is the assertion of \citet[Lemma 3.11]{sz2016b}, the proof of which 
is surprisingly technical. Another useful consequence is that the derivative 
of the unit normal vector is bounded, see \citet[Lemma 3.10]{sz2016b}.

\subsection{Definition of the transform and main results}

Our choice of the transform $G$ is 
\[
 G(x)=x+\alpha(p(x))\tilde \phi(x) \,,
\]
where
\[
\tilde \phi(x)=(x-p(x))\cdot
n(p(x))\|x-p(x)\|\,\phi\left(\frac{\|x-p(x)\|}{c}\right)\,.
\]
This should be compared to the 1D analog, equation \eqref{eq:transf1d}.
In \citet[Theorem 3.14 and Lemma 3.18]{sz2016b} it is proven that under the assumptions of Theorem \ref{thm:exun} below, $c$ can always be chosen sufficiently small, so that $G$ has a global inverse by Hadamard's global inverse function theorem \cite[Theorem 2.2]{ruzhansky2015}.
In 1D the constant $\alpha$ had the purpose of making sure that the jump of $\mu$ is offset by the jump of $G''$. In general dimension, $\alpha$ is defined on the hypersurface $\hypsurf$:
 \begin{align}\label{eq:alpha}
\alpha(\xi)=\lim_{h\to 0}\frac{\mu(\xi-h n(\xi))-\mu(\xi+hn(\xi))}{2 \|\sigma(\xi)^\top n(\xi)\|^2}\,,\qquad\xi\in \Theta\,.
\end{align}
Although $\alpha$ depends on the choice
of the normal unit vector, it is readily checked that $G$ does not.

We will need to make additional assumptions on $\mu$ and $\sigma$ to 
guarantee existence and
sufficient regularity of $\alpha$ and, a fortiori, of $G$.

It remains to show that It\^o's formula holds for $G^{-1}$. This follows from the following special case of \cite[Theorem 2.1]{peskir2007}.
\begin{theorem}[It\^o's formula]
Let $X$ be a $d$-dimensional It\^o process and let
$b:\R^{d-1}\longrightarrow\R$ be a $C^2$-function.
Let furthermore $f_1,f_2:\R^d\longrightarrow \R$ be $C^2$-functions such that 
the function $f:\R^d\longrightarrow \R$ defined by 
$$f(x)=f_1(x) \1_{x_{d}\le b(x_1,\ldots,x_{d-1})}+f_2(x) \1_{x_{d}> b(x_1,\ldots,x_{d-1})}$$ is in $C^1$.
Then It\^o's formula holds for $X$ and $f$.
\end{theorem}

We have the following existence and uniqueness result.

\goodbreak

\begin{theorem}[\citet{sz2016b}]\label{thm:exun}
Let the following assumptions hold:
\begin{itemize}\itemsep-\parsep
\item $\mu:\R^d\longrightarrow\R^d$ is piecewise Lipschitz with exceptional set $\Theta$;
\item $\Theta \in C^4$, has positive reach, and $n'', n'''$ are bounded;
\item $\sigma:\R^d \longrightarrow \R^{d\times d}$ is Lipschitz and $\|\sigma(\xi)^\top n(\xi)\|^2\ge c_0>0$ for all $\xi \in \Theta$;
\item $\mu,\sigma$ are bounded on $\Theta^\varepsilon$ for some $\varepsilon>0$;
\item $\mu,\sigma$ are such that $\alpha$, as described in \eqref{eq:alpha},
is well-defined and has bounded derivatives up to order 3.
\end{itemize}

 Then there exists a unique strong solution to \eqref{eq:SDE}.
\end{theorem}

We remark that the assumptions of Theorem \ref{thm:exun}
impose extra regularity on $\mu,\sigma$ only close to, and on $\Theta$. 
Away from $\Theta$ we basically have the classical Lipschitz requirements.
In analogy to the one-dimensional result, we have the following:
\begin{theorem}[\citet{sz2016b}]
Let the assumptions of Theorem \ref{thm:exun} hold.
Then, also in the multidimensional setting, Algorithm \ref{alg:GM} converges with strong order $1/2$.
\end{theorem}

\subsection{Example}

We apply our Algorithm \ref{alg:GM} to solve an example of an SDE where
the drift is discontinuous on the unit circle in the $\R^2$, i.e.~the exceptional set $\Theta=\{x \in \R^2 | x_1^2+x_2^2=1\}$, and the diffusion coefficient is degenerate.
Let
\begin{align*}
dX_t= \mu(X^1_t,X^2_t) dt + \sigma(X^1_t,X^2_t) dW_t\,,
\end{align*}
where
\begin{align*}
 \mu(x_1,x_2)&=\begin{cases}
(1,1)^\top, & x_1^2+x_2^2>1\\
(-x_1,x_2)^\top, & x_1^2+x_2^2\le1\,,
          \end{cases}
\\[1em]
 \sigma(x_1,x_2)&=\frac{1}{1+x_1^2+x_2^2} \left(
  \begin{array}{cc}
   x_1 & 0  \\
   x_2 & 0
    \end{array} \right)\,.
    \end{align*}

\begin{figure}[t]
	\centering
	\includegraphics[width=1.0\textwidth]{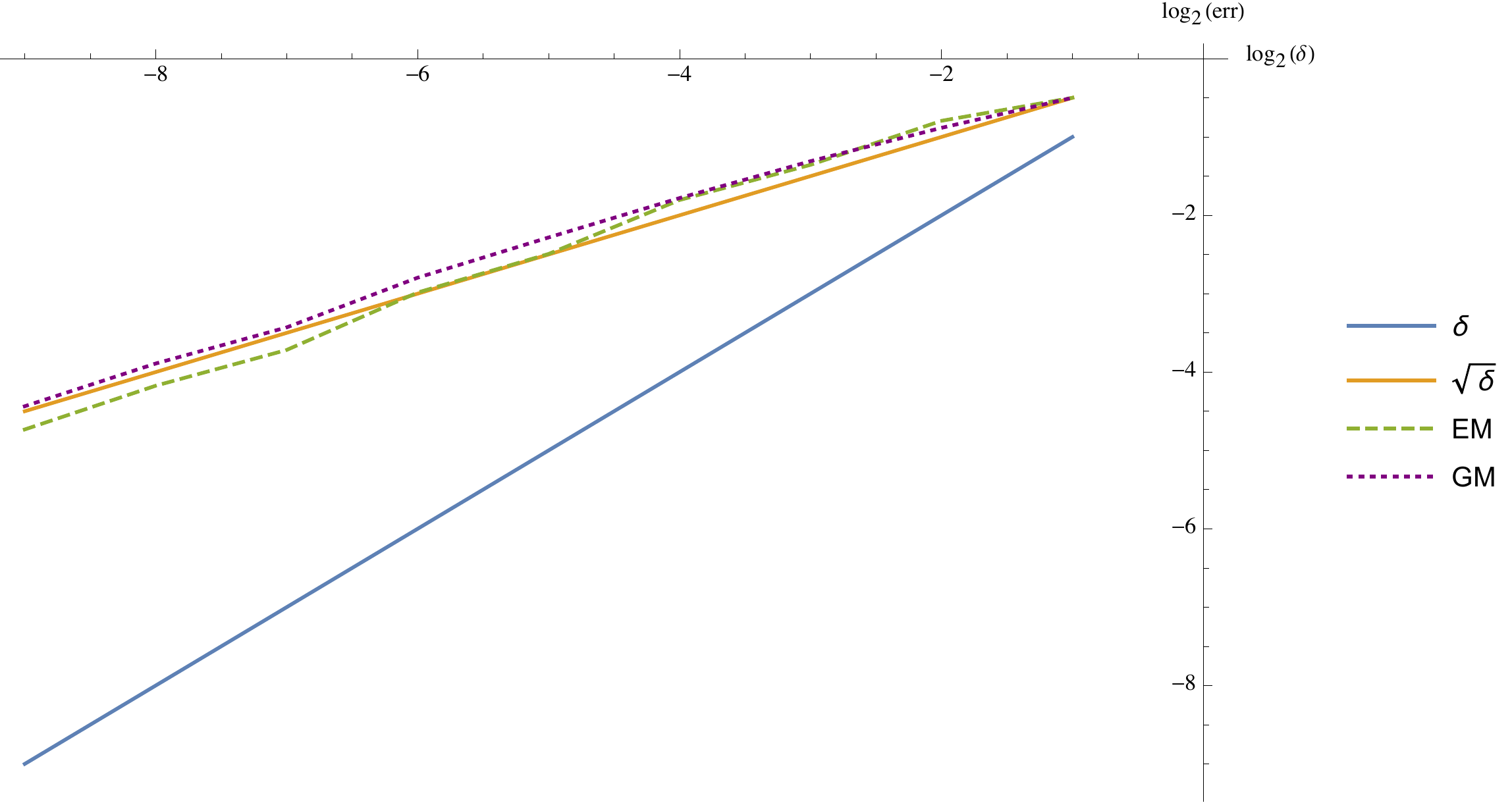}
	\caption{Estimated $L^2$-error of GM and EM.\label{fig:errorkreis}}
	\end{figure}

Figure \ref{fig:errorkreis} shows the estimated $L^2$-error of GM for this example.
We observe that GM shows the convergence behaviour we expect from our theoretical result, namely it converges as fast as $\delta^{1/2}$, i.e.~the purple dotted line has the same slope as the yellow line.
So in principle we could be satisfied.
We have constructed the first numerical method that is proven to converge for a rather general class of SDEs with discontinuous drift and we have established its convergence speed.
However, GM has two shortcomings.
First, it needs the geometrical structure of the set of discontinuities of the drift as an input. However, if for example the discontinuity stems from a discontinuous control policy in a stochastic optimal control problem, then this geometric structure for the optimal control is not explicitly known.
Finding the discontinuity of a function numerically is a problem of high complexity on its own.
Second, our method requires inversion of $G$ in each step. In 1D the inverse can be calculated explicitly, see Theorem \ref{thm:inverseexp}, but in general dimension, we have to resort to numerical inversion, which makes the calculation of a single path rather costly.

However, Figure \ref{fig:errorkreis} tells us even more. We observe that the green dashed line, which corresponds to the convergence speed of the EM method applied to our example, also has roughly the same slope as the yellow line.
This means that for our example and our range of $\delta$, the EM method seems to converge, too.
To deal with the issues raised above, it would be desirable to prove a positive strong convergence rate for the EM method.
This is what we are going to study in the next section.

\section{Convergence of the EM method}

We seek to estimate the mean square error of the EM approximation by considering the difference between GM and EM. Here, we only sketch the idea of the proof.

Let $\appx$ be the EM approximation of $X$.
Using that $X=G(Z)$, that $G^{-1}$ is Lipschitz, and that $(a+b)^2\le 2a^2+2b^2$, we estimate the mean square error of the EM approximation:
\begin{align*}
	& \E \Big(\sup_{0\le t\le T} \|X_t-\appx_t \|^2 \Big)
	 = \E \Big(\sup_{0\le t\le T} \|G^{-1}(Z_t)-G^{-1}(G (\appx_t )) \|^2 \Big)\\
	&\qquad \le 2L_{G^{-1}} ^2\E \Big(\sup_{0\le t\le T} \|Z_t- \appz_t \|^2 \Big)
+2L_{G^{-1}} ^2 \E \Big(\sup_{0\le t\le T} \|\appz_t-G (\appx_t ) \|^2 \Big)\,.
\end{align*}
With this we have decomposed the error into two error terms. The first term is the mean square error of the EM approximation of the solution to the transformed SDE. Since the transformed SDE has Lipschitz coefficients, 
the EM method converges with strong order $1/2$, i.e.
\[
\E\Big(\sup_{0\le t\le T}\|Z_t- \appz_t\|^2\Big)\le C \delta\,.
\]
For estimating
\begin{align*}
 \E \Big(\sup_{0\le t\le T} \|\appz_t-G (\appx_t ) \|^2 \Big)
 \end{align*}
the crucial estimate is the one of the drift. For this the main tasks are:

\begin{itemize}\itemsep-\parsep
 \item estimating the probability of the event $\Omega_\varepsilon$ that during one step the distance between the interpolation of the EM method and the previous EM step becomes greater than some given $\varepsilon>0$.
 Lemma 3.3 in \cite{sz2017c} states that
 \begin{align*}
  \P(\Omega_\varepsilon) \le C \exp\left(-\frac{\varepsilon}{\|\sigma\|_\infty \delta^{1/2}}\right)\,;
 \end{align*}
 \item estimating the occupation time of the Euler-Maruyama approximation of $X$
close to the hypersurface $\Theta$ by constructing a 1D process $Y$ that has
the same occupation time close to $0$ as $\appx$ has close to $\Theta$.
The process $Y$ is essentially a signed distance of $\appx$ from $\Theta$.
Again we make extensive use of the positive reach property of $\Theta$, which
guarantees regularity of a distance function.
Theorem 2.7 in \cite{sz2017c} says that
\begin{align*}
\int_0^{T} \P \left( \{\appx_s \in \Theta^\varepsilon \} \right) ds\le C\varepsilon\,.
\end{align*}
\end{itemize}
We are free to choose $\varepsilon$ as a function of the step-size $\delta$, and if we do so in an optimal way, we obtain the following convergence rate.

\begin{theorem}[\citet{sz2017c}]
Let the assumptions of Theorem \ref{thm:exun} hold, and let $\mu,\sigma$ be bounded.
 
 Then the Euler-Maruyama method converges with strong order $1/4-\othereps$ for arbitrarily small $\othereps>0$ to the solution of SDE \eqref{eq:SDE}.
\end{theorem}

Now the question arises why one would apply EM instead of GM, since GM has a much higher convergence speed. However, as already mentioned at the end of Section \ref{sec:ddim}, the computation of a single path with GM can be so slow, that obtaining comparable errors with GM can take more time for practical purposes. We refer to \cite{sz2017c} for more details.

\begin{figure}[t]
	\centering
	\includegraphics[width=1.0\textwidth]{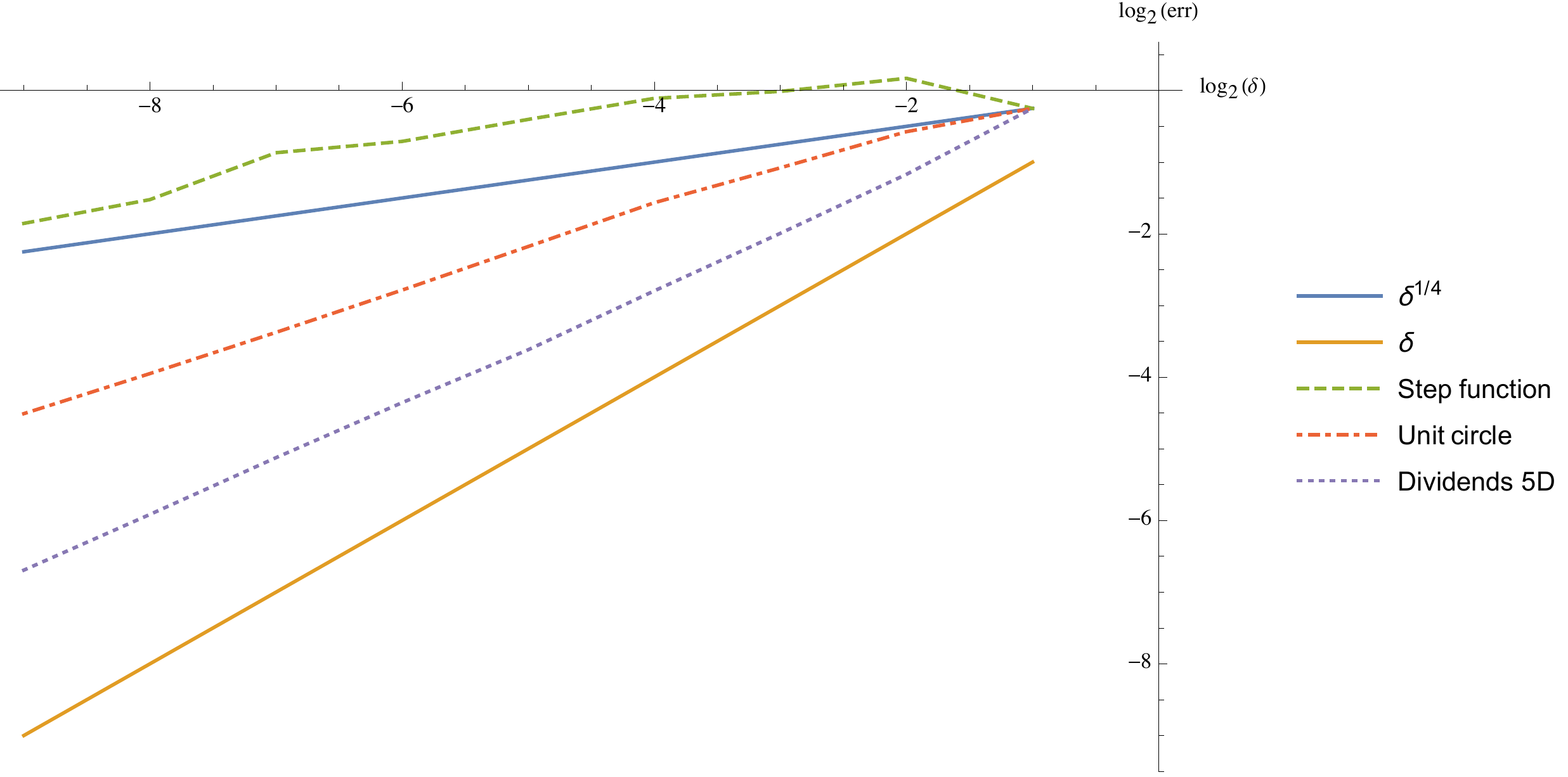}
	\caption{Estimated $L^2$-error of the EM approximation.\label{fig:errs}}
	\end{figure}

Figure \ref{fig:errs} shows the estimated $L^2$-error of the EM approximation for three examples: one where the drift is a certain step-function, a five-dimensional example from insurance mathematics (Dividends 5D), and the example from above where the drift is discontinuous on the unit circle.
We see that for the step-function example the convergence seems to be approximately as fast as $\delta^{1/4}$ for larger $\delta$, but for smaller $\delta$ the slope of the dashed green line seems to become steeper.
For the other two examples the EM method clearly converges at a higher rate for this example. This supports the claim from above that in many examples the EM method is the preferred choice.

\section*{Acknowledgements}

G.~Leobacher is supported by the Austrian Science Fund (FWF): Project F5508-N26, which is part of the Special Research Program ``Quasi-Monte Carlo Methods: Theory and Applications''.

M.~Sz\"olgyenyi is supported by the AXA Research Grant ``Numerical Methods for Stochastic Differential Equations with Irregular Coefficients with Applications in Risk Theory and Mathematical Finance".
This article was written while M.~Sz\"olgyenyi was affiliated with the Institute of Statistics and Mathematics, Vienna University of Economics and Business, Welthandelsplatz 1, 1020 Vienna, Austria, and supported by the Vienna Science and Technology Fund (WWTF): Project MA14-031.

\vspace{2em}

 \noindent G. Leobacher \Letter\\
 Institute for Mathematics and Scientific Computing, University of Graz, Heinrichstra\ss{}e 36, 8010 Graz, Austria\\
 gunther.leobacher@uni-graz.at\\

\noindent M. Sz\"olgyenyi \\
Seminar for Applied Mathematics, ETH Z\"urich, R\"amistra\ss{}e 101, 8092 Z\"urich, Switzerland\\
michaela.szoelgyenyi@sam.math.ethz.ch

\end{document}